# REDUCING VARIANCE IN UNIVARIATE SMOOTHING

By Ming-Yen Cheng, Liang Peng[1] and Jyh-Shyang Wu

*National Taiwan University, Georgia Institute of Technology and Tamkang University*

A variance reduction technique in nonparametric smoothing is proposed: at each point of estimation, form a linear combination of a preliminary estimator evaluated at nearby points with the coefficients specified so that the asymptotic bias remains unchanged. The nearby points are chosen to maximize the variance reduction. We study in detail the case of univariate local linear regression. While the new estimator retains many advantages of the local linear estimator, it has appealing asymptotic relative efficiencies. Bandwidth selection rules are available by a simple constant factor adjustment of those for local linear estimation. A simulation study indicates that the finite sample relative efficiency often matches the asymptotic relative efficiency for moderate sample sizes. This technique is very general and has a wide range of applications.

**1. Introduction.** Local linear modeling for nonparametric regression has many advantages and has become very popular. Hastie and Loader [16], Wand and Jones [27], Fan and Gijbels [11] and others have investigated extensively its theoretical and practical properties. To reduce the variance, for each point $x$, take a special linear combination of this local linear estimators at three equally spaced points around $x$. The linear combination satisfies certain moment conditions so that the asymptotic bias is unchanged. Below are a few specific features of this new estimator. First, both local and global automatic bandwidths can be easily obtained from those for the standard local linear estimator. Second, the asymptotic mean squared error is improved considerably and the amount of reduction is uniform across different locations, regression functions, designs and error distributions. Third, evidenced by a simulation study, the reduction in asymptotic variance is effectively

Received May 2004; revised May 2006.
[1]Supported by NSF Grant DMS-04-03443.

projected to finite samples. Fourth, the estimators admit simple forms and only slightly increase the amount of computation time. Finally, many advantages of the local linear estimator, for example, design adaptivity and automatic boundary correction, are retained.

In kernel density estimation, Kogure [18] studied a related topic on polynomial interpolation and obtained some preliminary results. The purpose is to find optimal allocation of the interpolation points by minimizing the asymptotic mean integrated squared error. Higher-order polynomial interpolation does not change the asymptotic bias but alters the asymptotic variance in a nonhomogeneous way unless the spacings are all large enough, which was assumed in seeking the optimal allocation.

Fan [9] showed that the local linear estimator is minimax optimal among all linear estimators. The proposed estimators are linear and have smaller asymptotic mean squared errors than the local linear estimator. There is no conflict between these two results since the assumptions are different. In [9], the maximum risk of a linear estimator is taken over a class of regression functions that can be approximated well linearly. The asymptotic results for our estimators require the slightly more stringent condition that the regression function has a bounded, continuous second derivative at the point of estimation.

There is an extensive literature on modifications and improvements of kernel and local linear estimators. Many of them are aimed at bias reduction, for example, Abramson [1], Samiuddin and El-Sayyad [22], Jones, Linton and Nielsen [17] and Choi and Hall [7]. Choi and Hall [7] exploit the idea of taking linear combinations of some local linear estimators as well. The key difference is that we take linear combinations to maintain bias but reduce variance while Choi and Hall [7] take linear combinations to reduce bias while essentially maintaining variance. More specifically, (i) they use local linear estimates at points symmetrically located around $x$ and we do not, (ii) they take a convex combination and we employ different constraints on the coefficients, (iii) bandwidth selection is straightforward in our case but is much more complicated for theirs, and (iv) our method applies to local constant fitting estimators but theirs does not. There are very few variance reduction techniques available. An error-dependent technique in local linear smoothing was suggested by Cheng and Hall [3]. The amount of variance reduction depends on the error distribution. Cheng and Hall [4] introduced an adaptive line integral of the bivariate kernel density estimate to reduce variance. It does not apply to the univariate case or regression estimation. These two methods require explicitly or implicitly a pilot estimation of the function or its derivatives. Our procedure is based on a completely different idea and effectively reduces variance without pilot estimation.

Section 2 presents the methodology and Section 3 provides results on the asymptotic mean squared error and coverage accuracy. Section 4 discusses some practical issues, including bandwidth selection and implementation. Section 5 contains a numerical study. Possible generalizations and applications of the variance reduction technique are discussed in Section 6. Proofs are given in Section 7.

## 2. Methodology.

2.1. *Local linear regression.* Suppose that independent bivariate observations $(X_1, Y_1), \ldots, (X_n, Y_n)$ are sampled from the regression model

$$Y = m(X) + \sigma(X)\varepsilon,$$

where $\sigma(X)$ is the conditional variance and the random error $\varepsilon$, independent of $X$, has zero mean and unit variance. Kernel estimators of $m(x) = E(Y|X=x)$ include the Nadaraya–Watson, Gasser–Müller and local polynomial estimators; see [8, 25, 27]. Let $K$ be a kernel function and $h > 0$ be a bandwidth. The local linear regressor is defined as

$$(2.1) \qquad \frac{S_{n,2}(x)T_{n,0}(x) - S_{n,1}(x)T_{n,1}(x)}{S_{n,0}(x)S_{n,2}(x) - S_{n,1}(x)S_{n,1}(x)},$$

where $S_{n,l}(x) = h\sum_{i=1}^{n}(x - X_i)^l K_h(x - X_i)$, $l = 0, 1, 2$, $T_{n,l}(x) = h\sum_{i=1}^{n}(x - X_i)^l K_h(x - X_i)Y_i$, $l = 0, 1$, and $K_h(t) = K(t/h)/h$. This estimator is obtained by solving a local linear least squares problem. Its appealing theoretical and numerical properties were discussed by Fan [9] and Fan and Gijbels [11], among others. Further, it has become very popular, widely used in applications and implemented in statistical software. When the denominator is close to zero, it exhibits a rather unstable numerical behavior. This is particularly problematic for small samples or sparse designs. Remedies include modifications proposed by Seifert and Gasser [23], Cheng, Hall and Titterington [5], Hall and Turlach [15] and Seifert and Gasser [24]. Although we focus on the local linear estimator denoted by $\widehat{m}(x)$, our variance reduction methods given below can be applied to any such modification.

2.2. *Variance reduced estimators.* A motivation of our variance reduction strategy is to incorporate more data points in the regression estimation in such a way that the first-order term in the asymptotic bias remains unchanged. To do this, at each $x$ we construct a linear combination of the local linear estimators at some equally spaced points near $x$, with the linear coefficients satisfying certain moment conditions derived from the asymptotic bias expansions. Then variance reduction is achieved since the three preliminary estimators are correlated and the correlation coefficients are less than 1. In this context, we fix the number of nearby points at three since that is the minimal requirement specified by the moment conditions and using more than three yields complex solutions. In addition, the grid points are taken so that the final estimator is the simplest and most efficient in the mean squared error sense.

Formally, for any $x$, define three equally spaced points $\alpha_{x,j} = x - (r+1-j)\omega_n$, $j = 0, 1, 2$, where $r \in (-1, 1)$ and $\omega_n = \delta h$ for some constant $\delta > 0$. The shift parameter $r$ determines the location of the leftmost point $\alpha_{x,0}$ relative to $x$ and $\omega_n$ is the spacing of the grid. Construct a quadratic interpolation of the local linear estimators $\widehat{m}(\alpha_{x,0}), \widehat{m}(\alpha_{x,1})$ and $\widehat{m}(\alpha_{x,2})$ and then estimate $m(x)$ by the value of the interpolated curve at $x$:

$$\widetilde{m}_q(x) = \sum_{j=0,1,2} A_j(r)\widehat{m}(x - (r+1-j)\omega_n),$$

where the coefficients depend only on $r$ and are given by

$$A_0(r) = r(r-1)/2, \qquad A_1(r) = (1-r^2), \qquad A_2(r) = r(r+1)/2.$$

Then the moment conditions (7.4) are satisfied so that $\widetilde{m}_q(x)$ has the same asymptotic bias as $\widehat{m}(x)$. Theorem 1 and Proposition 1 show that $\widetilde{m}_q(x)$ has a smaller asymptotic variance than $\widehat{m}(x)$ and the asymptotic variances differ from each other by a constant factor depending on the bin-width parameter $\delta$, the shift parameter $r$ and the kernel $K$. Given $K$ and $\delta$, the optimal values of $r$ that minimize this constant factor are $r = \pm\sqrt{1/2}$, which give

$$\widetilde{m}_\pm(x) = \sum_{j=0,1,2} A_j(\pm 1/\sqrt{2})\widehat{m}(x - (\pm 1/\sqrt{2} + 1 - j)\omega_n).$$

Since $\widetilde{m}_+(x)$ uses more data information on the left-hand side of $x$ than on the right-hand side, the curve estimate $\widetilde{m}_+(\cdot)$ tends to shift to the right of $m(\cdot)$. Similarly, $\widetilde{m}_-(\cdot)$ tends to shift to the left of $\widehat{m}(\cdot)$. This symmetry suggests taking the average

$$\widetilde{m}_a(x) = \{\widetilde{m}_+(x) + \widetilde{m}_-(x)\}/2$$

as our final estimator. Sections 3 and 5 demonstrate that, compared to $\widetilde{m}_\pm(x)$, $\widetilde{m}_a(x)$ further improves the asymptotic and finite sample efficiencies.

When $\text{Supp}(m)$ is bounded, $\text{Supp}(m) = [0, 1]$ say, since each of $\widetilde{m}_\pm(x)$ and $\widetilde{m}_a(x)$ requires values of $\widehat{m}$ at points around $x$, a $\delta$ value,

(2.2) $\qquad \delta(x) = \min\{\delta, x/(\sqrt{1/2}+1)h, (1-x)/(\sqrt{1/2}+1)h\},$

which depends on the distances from $x$ to the boundary points 0 and 1, is used so that $\widetilde{m}_a(x)$ is defined for every $x \in [0, 1]$.

Our estimators $\widetilde{m}_\pm(x)$ and $\widetilde{m}_a(x)$ are simple linear combinations of local linear estimators evaluated at nearby points. No pilot estimation is involved in this variance reduction. Therefore, in finite samples, the asymptotic efficiencies are achieved at relatively small $n$. For the same reasons, they share many advantages enjoyed by the local linear estimator, for example automatic boundary correction and design adaptivity; see [11]. It is shown in Section 3 that each of $\widetilde{m}_\pm(x)$ and $\widetilde{m}_a(x)$ reduces the asymptotic variance uniformly across different regression functions, different designs and different error distributions. Then procedures such as design planning employed in applications of $\widehat{m}$ apply to our estimators in general. Thus the proposed estimators are rather user friendly.

2.3. *Confidence intervals.* Consider constructing confidence intervals for $m(x)$ based on $\widehat{m}(x)$ and $\widetilde{m}_q(x)$. Let $\nu_{ij} = \int s^i K(s)^j \, ds$ and $\widetilde{\nu}_{0l} = \int \{\sum_{i=0,1,2} A_i(r) K(s+i\delta)\}^l \, ds$. Define $w_{ijk}(x) = n^{-1} h^{j-i-1} \sum_{l=1}^n (x-X_l)^i \{K_h(x-X_l)\}^j \times \{Y_l - m(x)\}^k$ and $\widehat{\sigma}^2(x) = n^{-1} \sum_{l=1}^n K_h(x-X_l)\{Y_l - \widehat{m}(x)\}^2 / w_{010}(x)$. Then both

$$T_n(x) = \frac{\sqrt{nh}}{\sqrt{\nu_{02}}} \frac{\widehat{m}(x) - m(x)}{\sqrt{\widehat{\sigma}^2(x)/w_{010}(x)}}, \qquad T_n^*(x) = \frac{\sqrt{nh}}{\sqrt{\widetilde{\nu}_{02}}} \frac{\widetilde{m}_q(x) - m(x)}{\sqrt{\widehat{\sigma}^2(x)/w_{010}(x)}}$$

are asymptotically $N(0,1)$ distributed when $nh^5 \to 0$, and we have the following one-sided confidence intervals for $m(x)$:

(2.3) $\quad I_\beta = (\widehat{m}(x) - z_\beta \{\widehat{\sigma}^2(x)/w_{010}(x)\}^{1/2} \nu_{02}^{1/2} (nh)^{-1/2}, \infty),$

(2.4) $\quad \widetilde{I}_\beta = (\widetilde{m}_q(x) - z_\beta \{\widehat{\sigma}^2(x)/w_{010}(x)\}^{1/2} \widetilde{\nu}_{02}^{1/2} (nh)^{-1/2}, \infty),$

where $z_\beta$ satisfies $P\{N(0,1) \leq z_\beta\} = \beta$.

**3. Theoretical results.** The following conditions are needed for asymptotic analysis of the estimators.

ASSUMPTION $(C_K)$. $K$ is a symmetric density function with compact support.

ASSUMPTION $(C_{m,f})$.

1. $m''(\cdot)$, the second derivative of $m(\cdot)$, is bounded and continuous at $x$;
2. the density function $f(\cdot)$ of $X$ satisfies $f(x) > 0$ and $|f(x) - f(y)| \leq c|x-y|^\alpha$ for some $0 < \alpha < 1$;
3. $\sigma^2(\cdot)$ is bounded and continuous at $x$.

3.1. *Asymptotic mean squared error.* Fan [9] showed that, under the conditions $(C_{m,f})$, $(C_K)$, $h \to 0$ and $nh \to \infty$ as $n \to \infty$,

(3.1) $\quad \widehat{m}(x) = \frac{S_{n,2}(x) T_{n,0}(x) - S_{n,1}(x) T_{n,1}(x)}{S_{n,0}(x) S_{n,2}(x) - S_{n,1}(x) S_{n,1}(x) + n^{-2}},$

a modification of (2.1) that admits asymptotic unconditional variance, has

$$\text{E}\{\widehat{m}(x)\} - m(x) = \tfrac{1}{2}m''(x)\nu_{20}h^2 + o\{h^2 + (nh)^{-1/2}\}, \tag{3.2}$$

$$\text{Var}\{\widehat{m}(x)\} = \frac{\sigma^2(x)}{nhf(x)}\nu_{02} + o\{h^4 + (nh)^{-1}\}. \tag{3.3}$$

THEOREM 1. *Suppose that $\delta > 0$ is a constant. Assume $(C_{m,f})$, $(C_K)$, $h \to 0$ and $nh \to \infty$ as $n \to \infty$. Then*

$$\text{E}\{\widetilde{m}_q(x)\} - m(x) = \tfrac{1}{2}m''(x)\nu_{20}h^2 + o\{h^2 + (nh)^{-1/2}\}, \tag{3.4}$$

$$\text{Var}\{\widetilde{m}_q(x)\} = \frac{\sigma^2(x)}{nhf(x)}\widetilde{\nu}_{02} + o\{h^4 + (nh)^{-1}\}. \tag{3.5}$$

Note that $\widetilde{\nu}_{02} = \nu_{02} - r^2(1-r^2)C(\delta)$, where

$$C(\delta) = 1.5C(0,\delta) - 2C(0.5,\delta) + 0.5C(1,\delta)$$

with $C(a,\delta) = \int K(t-a\delta)K(t+a\delta)\,dt$.

PROPOSITION 1. *The quantity $C(\delta)$ has the following properties:*

(a) *For any symmetric kernel function $K$, $C(\delta) \geq 0$ for any $\delta \geq 0$.*
(b) *If $K$ has a unique maximum and is concave, then $C(\delta)$ is increasing in $\delta > 0$.*

REMARK 1. From (3.2) and (3.4), $\widetilde{m}_q(x)$ and $\widehat{m}(x)$ have the same asymptotic bias. From (3.3), (3.5) and Proposition 1, the asymptotic variance of $\widetilde{m}_q(x)$ is smaller than that of $\widehat{m}(x)$ by the amount $\{nhf(x)\}^{-1}\sigma^2(x)r^2(1-r^2)C(\delta)$. Note that $0 < r^2(1-r^2) \leq 1/4$ for any $r \in (-1,1) \setminus \{0\}$ and attains its maximum $1/4$ at $r = \pm\sqrt{1/2}$. Therefore, for any $\delta > 0$, the optimal choices of $r$ are $r = \pm\sqrt{1/2}$, which yield $\widetilde{m}_\pm(x)$.

REMARK 2. A generalization of $\widehat{m}_q(x)$, based on local linear estimators at $z_0 = x - r\delta h$, $z_1 = x - (r-k)\delta h$ and $z_2 = x - (r-k-1)\delta h$, for some $0 < k \neq 1$, is $\sum_{j=0,1,2} B_j(r)\widehat{m}(z_j)$, where $B_0(r) = r(r-1)/k(k+1)$, $B_1(r) = -(r+k)(r-1)/k$ and $B_2(r) = r(r+k)/(k+1)$. It has the same asymptotic bias as $\widehat{m}(x)$ and asymptotic variance $\{nhf(x)\}^{-1}\sigma^2(x)\tau(\delta,r,k)$, where

$$\tau(\delta,r,k) = \nu_{02} \sum_{j=0,1,2} B_j(r)^2 + 2B_0(r)B_1(r)C(k,\delta/2)$$
$$+ 2B_0(r)B_2(r)C(k+1,\delta/2) + 2B_1(r)B_2(r)C(1,\delta/2).$$

In $\tau(\delta,r,k)$, $K$ interacts with $r$, $k$ and $\delta$ and there is no explicit value of $r$ minimizing $\tau(\delta,r,k)$ for given $\delta$ and $k$.

COROLLARY 1. *Under the conditions in Theorem 1, as $n \to \infty$,*

(3.6) $\quad \mathrm{E}\{\widetilde{m}_{\pm}(x)\} - m(x) = \frac{1}{2}m''(x)\nu_{20}h^2 + o\{h^2 + (nh)^{-1/2}\},$

(3.7) $\quad \mathrm{Var}\{\widetilde{m}_{\pm}(x)\} = \frac{\sigma^2(x)}{nhf(x)}\left\{\nu_{02} - \frac{C(\delta)}{4}\right\} + o\{h^4 + (nh)^{-1}\}.$

Theorem 2 can be proved using arguments similar to the proof of Theorem 1.

THEOREM 2. *Under the conditions in Theorem 1, as $n \to \infty$,*

(3.8) $\mathrm{E}\{\widetilde{m}_a(x)\} - m(x) = \frac{1}{2}m''(x)\nu_{20}h^2 + o\{h^2 + (nh)^{-1/2}\},$

(3.9) $\quad \mathrm{Var}\{\widetilde{m}_a(x)\} = \frac{\sigma^2(x)}{nhf(x)}\left\{\nu_{02} - \frac{C(\delta)}{4} - \frac{D(\delta)}{2}\right\} + o\{h^4 + (nh)^{-1}\},$

*where*

$$\begin{aligned} D(\delta) = {} & \nu_{02} - \tfrac{1}{4}C(\delta) \\ & - \tfrac{1}{16}\{4(1+\sqrt{2})C(\sqrt{2}-1,\delta/2) \\ & \quad + (3+2\sqrt{2})C(2-\sqrt{2},\delta/2) \\ & \quad + 2C(\sqrt{2},\delta/2) + 4(1-\sqrt{2})C(\sqrt{2}+1,\delta/2) \\ & \quad + (3-2\sqrt{2})C(\sqrt{2}+2,\delta/2)\}. \end{aligned}$$

PROPOSITION 2. *The quantity $D(\delta)$ in (3.9) is nonnegative for any $\delta \geq 0$.*

REMARK 3. If $m^{(4)}(x)$ exists, then the second-order term in the bias of $\widetilde{m}_a(x)$ is $O(h^4)$, smaller than those of $\widehat{m}(x)$ and $\widehat{m}_{\pm}(x)$.

REMARK 4. Suppose that $\mathrm{supp}(K) = [-1,1]$. Then, for $\delta \geq 2$, $C(\delta) = (3/2)\nu_{02}$ and

(3.10) $$\begin{aligned} \mathrm{Var}\{\widetilde{m}_q(x)\} &\approx \{1 - \tfrac{3}{2}r^2(1-r^2)\}\mathrm{Var}\{\widehat{m}(x)\}, \\ \mathrm{Var}\{\widetilde{m}_{\pm}(x)\} &\approx \tfrac{5}{8}\mathrm{Var}\{\widehat{m}(x)\}. \end{aligned}$$

For $\delta \geq 2/(\sqrt{2}-1)$, $D(\delta) = (5/8)\nu_{02}$ and

(3.11) $\quad\quad\quad\quad\quad\quad \mathrm{Var}\{\widetilde{m}_a(x)\} \approx \tfrac{5}{16}\mathrm{Var}\{\widehat{m}(x)\}.$

If $K$ is infinitely supported, then (3.10) and (3.11) hold for sufficiently large $\delta$.

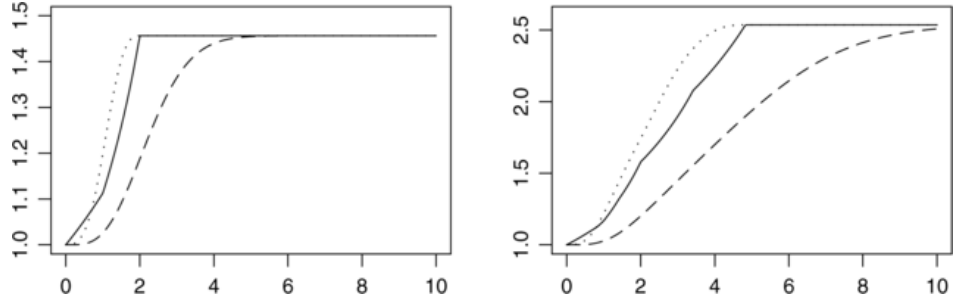

FIG. 1. *AMSE relative efficiency. The left and right panels respectively plot $\gamma_q(\delta)$ and $\gamma_a(\delta)$ against $\delta$ for the Uniform (solid), Epanechnikov (dotted) and Normal (dashed) kernels.*

REMARK 5. If $\delta = o(1)$, then $\widetilde{m}_q(x)$, $\widetilde{m}_\pm(x)$ and $\widetilde{m}_a(x)$ all have the same asymptotic bias and variance as $\widehat{m}(x)$. If $\delta \to \infty$ with $\delta = o(h^{-1/3})$ and $m'''(x)$ exists, then the biases of $\widetilde{m}_q(x)$ and $\widetilde{m}_\pm(x)$ remain the same as in (3.4) and (3.6) and the variances are as in (3.10). If $\delta \to \infty$ with $\delta = o(h^{-1/2})$ and $m^{(4)}(x)$ exists, then the bias and variance of $\widetilde{m}_a(x)$ are as in (3.8) and (3.11).

From Corollary 1, the pointwise (global) asymptotic efficiency, in terms of asymptotically optimal (integrated) MSE, of $\widetilde{m}_\pm(x)$ relative to $\widehat{m}(x)$ is

$$(3.12) \qquad \gamma_q(\delta) = \{\nu_{02} - C(\delta)/4\}^{-4/5} \nu_{02}^{4/5}.$$

In addition, Theorem 2 implies that the pointwise or global asymptotic relative efficiency of $\widetilde{m}_a(x)$ compared to $\widehat{m}(x)$ is

$$(3.13) \qquad \gamma_a(\delta) = \{\nu_{02} - C(\delta)/4 - D(\delta)/2\}^{-4/5} \nu_{02}^{4/5}.$$

Both $\gamma_q(\delta)$ and $\gamma_a(\delta)$ depend only on $K$ and $\delta$ and do not involve the regression function $m$, the design density $f$ or the error distribution. Thus, asymptotically, the variance reduction methods do not interfere with these factors. Figure 1 plots $\gamma_q(\delta)$ and $\gamma_a(\delta)$ against $\delta$ when $K$ is the Uniform, $K(u) = 0.5 I(|u|<1)$, Epanechnikov, $K(u) = 0.75(1-u^2)I(|u|<1)$ or Normal, $K(u) = (2\pi)^{-1/2} \exp(-u^2/2)$, kernel.

3.2. *Coverage probability.* Coverage probabilities of the confidence intervals $I_\beta$ and $\widetilde{I}_\beta$ for $m(x)$, given in (2.3) and (2.4) and constructed based on $\widehat{m}(x)$ and $\widetilde{m}_q(x)$, are analyzed and compared as follows.

THEOREM 3. *Assume $(C_K)$, $(C_{m,f})$, $h = o(n^{-1/5})$ and $nh/\log n \to \infty$ as $n \to \infty$. Then*

$$P\{m(x) \in I_\beta\}$$

$$
\begin{aligned}
(3.14) \quad &= \beta + (nh^5)^{1/2} 4^{-1} \nu_{21} \nu_{02}^{-1/2} \sigma^{-1}(x) \\
&\quad \times f^{1/2}(x) m''(x) (z_\beta^2 - 3) \phi(z_\beta) \\
&\quad - (nh)^{-1/2} 6^{-1} \nu_{02}^{-3/2} \sigma^{-3}(x) f^{-1/2}(x) V_3(x) \\
&\quad \times \{\nu_{03}(z_\beta^2 - 1) - 3\nu_{02}^2 z_\beta^2\} \phi(z_\beta) \\
&\quad + O\{(1/\sqrt{nh} + h)^2\},
\end{aligned}
$$

$$
\begin{aligned}
&P\{m(x) \in \widetilde{I}_\beta\} \\
(3.15) \quad &= \beta + (nh^5)^{1/2} 4^{-1} \nu_{21} \widetilde{\nu}_{02}^{-1/2} \sigma^{-1}(x) \\
&\quad \times f^{1/2}(x) m''(x) (z_\beta^2 - 3) \phi(z_\beta) \\
&\quad - (nh)^{-1/2} 6^{-1} \widetilde{\nu}_{02}^{-3/2} \sigma^{-3}(x) f^{-1/2}(x) V_3(x) \\
&\quad \times \{\widetilde{\nu}_{03}(z_\beta^2 - 1) - 3\widetilde{\nu}_{02}^2 z_\beta^2\} \phi(z_\beta) \\
&\quad + O\{(1/\sqrt{nh} + h)^2\},
\end{aligned}
$$

where $V_3(x) = E[\{Y - m(x)\}^3 | X = x]$.

COROLLARY 2. *Assume the conditions of Theorem 3. If $\{m''(x)(z_\beta^2 - 3)\}^{-1}\{\nu_{03}(z_\beta^2 - 1) - 3\nu_{02}^2 z_\beta^2\} < 0$ and $\{m''(x)(z_\beta^2 - 3)\}^{-1}\{\widetilde{\nu}_{03}(z_\beta^2 - 1) - 3\widetilde{\nu}_{02}^2 z_\beta^2\} < 0$, then*

$$
(3.16) \quad \begin{aligned}
\Gamma_\beta(\delta, r) &\equiv \lim_{n \to \infty} \frac{\min_h |P\{m(x) \in I_\beta\} - \beta|}{\min_h |P\{m(x) \in \widetilde{I}_\beta\} - \beta|} \\
&= \left\{\frac{\nu_{03}(z_\beta^2 - 1) - 3\nu_{02}^2 z_\beta^2}{\widetilde{\nu}_{03}(z_\beta^2 - 1) - 3\widetilde{\nu}_{02}^2 z_\beta^2}\right\}^{5/6} \left\{\frac{\widetilde{\nu}_{02}}{\nu_{02}}\right\}^{4/3}.
\end{aligned}
$$

Figure 2 plots $\Gamma_{0.95}(\delta, r)$ against $r$ and $\delta$ for the Uniform, Normal and Epanechnikov kernels. Similarly to the mean squared error comparison, $\Gamma_{0.95}(\delta, r)$ is always greater than or equal to 1, nondecreasing in $\delta$ both for $r = \pm 1/\sqrt{2}$ and for the best $r$, and settles at an upper limit. Also, $\Gamma_{0.95}(\delta, \pm 1/\sqrt{2})$ is very close to the optimal $\max_r \Gamma_{0.95}(\delta, r)$. Results for other $\beta$ values are similar. The limit $\lim_{\delta \to \infty} \max_r \Gamma_\beta(\delta, r)$ is roughly 1.15 if $\beta \geq 0.9$, 1.18 if $\beta = 0.85$, and 1.24 when $\beta = 0.75$ or 0.7. If $\beta = 0.75$ or 0.7, $\Gamma_\beta(\delta, r)$ can be less than 1, but such small confidence levels are rarely used in practice. For any fixed $\delta$, taking $r = \pm 1/\sqrt{2}$, as employed by $\widetilde{m}_\pm(x)$, the coverage accuracy is always not far from optimal. Table 1 gives $\Gamma_\beta(\delta, \pm 1/\sqrt{2})$ for different $\beta$ and $\delta$ values. There are significant gains in terms of coverage accuracy for $\delta \geq 1$ when using the Epanechnikov kernel.

10   M.-Y. CHENG, L. PENG AND J.-S. WU

TABLE 1
*Coverage accuracy ratio*

| $\delta$ | 0.6 | 0.8 | 1.0 | 1.2 | 1.6 | 2.0 |
|---|---|---|---|---|---|---|
| Uniform | 1.035 | 1.047 | 1.060 | 1.080 | 1.116 | 1.139 |
| | 1.031 | 1.042 | 1.054 | 1.074 | 1.115 | 1.152 |
| | 1.027 | 1.037 | 1.047 | 1.067 | 1.113 | 1.167 |
| | 1.022 | 1.031 | 1.039 | 1.060 | 1.112 | 1.184 |
| Epanechnikov | 1.024 | 1.045 | 1.067 | 1.088 | 1.123 | 1.136 |
| | 1.022 | 1.045 | 1.072 | 1.099 | 1.139 | 1.151 |
| | 1.021 | 1.045 | 1.078 | 1.110 | 1.156 | 1.167 |
| | 1.019 | 1.045 | 1.084 | 1.124 | 1.177 | 1.185 |
| Normal | 1.001 | 1.003 | 1.006 | 1.011 | 1.027 | 1.047 |
| | 1.001 | 1.004 | 1.008 | 1.015 | 1.035 | 1.059 |
| | 1.002 | 1.004 | 1.010 | 1.019 | 1.042 | 1.072 |
| | 1.002 | 1.006 | 1.013 | 1.023 | 1.051 | 1.086 |

$\Gamma_\beta(\delta, \pm 1/\sqrt{2})$ for the Uniform, Epanechnikov and Normal kernels. From the top, the rows respectively represent $\beta = 0.95$, 0.9, 0.85 and 0.8.

## 4. Implementation.

4.1. *Bandwidth selection.* Bandwidth choice is most crucial in kernel smoothing and dominates the performance. Compared to $\widehat{m}(x)$, $\widetilde{m}_\pm(x)$ and $\widetilde{m}_a(x)$ do not further complicate the bandwidth selection problem. This property does not necessarily hold for other variance or bias reducing modifications.

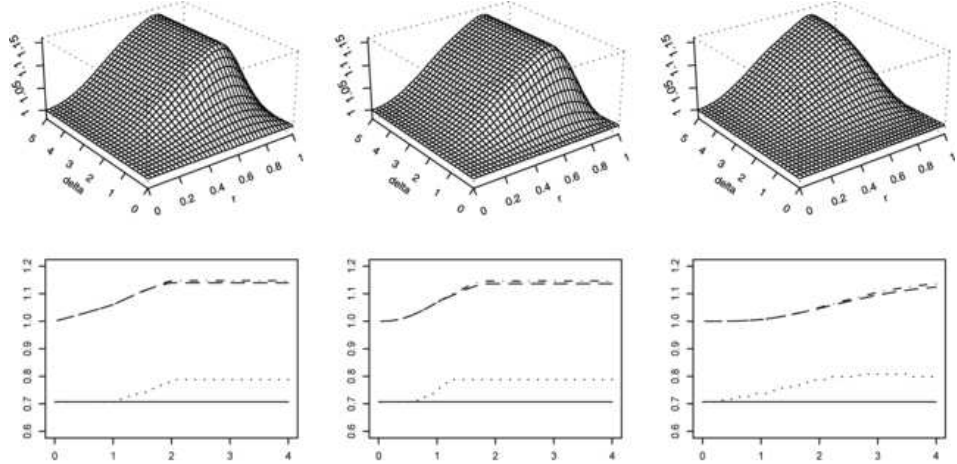

FIG. 2. *Coverage accuracy relative efficiency. Top row*: perspective plots of $\Gamma_{0.95}(\delta, r)$. *Bottom row*: plotted against $\delta$ are $\arg\max_r \Gamma_{0.95}(\delta, r)$ (*dotted*), $1/\sqrt{2}$ (*solid*), $\max_r \Gamma_{0.95}(\delta, r)$ (*dotted-dashed*) and $\Gamma_{0.95}(\delta, \pm 1/\sqrt{2})$ (*dashed*). *From left, the columns correspond to the Uniform, Epanechnikov and Normal kernels.*

First, consider local bandwidths which are needed when the underlying population has sharp features in the regression, design or error distribution. The optimal local bandwidth that minimizes the asymptotic mean squared error (AMSE) of $\widehat{m}(x)$ is

$$h_0(x) = \{\sigma^2(x)\nu_{02}\}^{1/5}\{nf(x)m''(x)^2\nu_{20}^2\}^{-1/5},$$

which gives the optimal AMSE

(4.1) $$1.25\{m''(x)^2\nu_{20}^2\sigma^8(x)/f^4(x)\}^{1/5}\nu_{02}^{4/5}n^{-4/5}.$$

For either of $\widetilde{m}_\pm(x)$, the optimal local bandwidth is

(4.2) $$h_1(x) = \{\nu_{02} - C(\delta)/4\}^{1/5}\nu_{02}^{-1/5}h_0(x),$$

yielding the optimal AMSE

(4.3) $$1.25\{m''(x)^2\nu_{20}^2\sigma^8(x)/f^4(x)\}^{1/5}\{\nu_{02} - C(\delta)/4\}^{4/5}n^{-4/5}.$$

The bandwidth that minimizes the AMSE of $\widetilde{m}_a(x)$ is

(4.4) $$h_a(x) = \{\nu_{02} - C(\delta)/4 - D(\delta)/2\}^{1/5}\nu_{02}^{-1/5}h_0(x),$$

giving the optimal AMSE

(4.5) $$1.25\{m''(x)^2\nu_{20}^2\sigma^8(x)/f^4(x)\}^{1/5}\{\nu_{02} - C(\delta)/4 - D(\delta)/2\}^{4/5}n^{-4/5}.$$

An implication of (4.2) and (4.4) is that, for any given $\delta$ and $K$, after adjusting by an appropriate constant multiplier, any local data-driven bandwidth designed for $\widehat{m}(x)$ is readily applicable to each of $\widetilde{m}_\pm(x)$ and $\widetilde{m}_a(x)$. In addition, regardless of what the regression, design or error distribution is, the multiplicative factors in (4.2) and (4.4) all remain the same for different $x$ values. Automatic local bandwidth selectors for local linear regression include those of [2, 10, 20].

Alternatively, for a generic kernel estimator $\bar{m}$ of $m$, consider the optimal global bandwidth that minimizes the asymptotic mean integrated squared error (AMISE), that is, the first-order term in the mean integrated squared error $\mathrm{E}\int\{\bar{m}(x) - m(x)\}^2 f(x)\,dx$. The optimal global bandwidths $h_0$, $h_1$ and $h_a$ respectively minimizing the AMISE's of $\widehat{m}$, $\widetilde{m}_\pm$ and $\widetilde{m}_a$ admit the relations

(4.6) $$h_1 = \{\nu_{02} - C(\delta)/4\}^{1/5}\nu_{02}^{-1/5}h_0,$$

(4.7) $$h_a = \{\nu_{02} - C(\delta)/4 - D(\delta)/2\}^{1/5}\nu_{02}^{-1/5}h_0.$$

Again, any automatic global bandwidth for the local linear estimator (see, e.g., [13, 21]) can be accordingly adjusted for implementation of $\widetilde{m}_\pm$ and $\widetilde{m}_a$.

Note that the constant factors in (4.2) and (4.6) are the same, and those in (4.4) and (4.7) are equal. Thus, bandwidth selection for our estimators is very simple. This advantage arises from the fact that $r$ is kept fixed for all $x$, thereby enhancing AMSE performance uniformly across different $x$, regressions, designs and error distributions.

4.2. *Bin width.* In the construction of $\widetilde{m}_{\pm}(x)$ and $\widetilde{m}_a(x)$, the equally spaced points $\alpha_{x,0}$, $\alpha_{x,1}$ and $\alpha_{x,2}$ have bin width $\delta > 0$. If $\delta = o(1)$, then there is no variance reduction; see Remark 5. If $\delta$ diverges with $\delta = o(h^{-1/3})$, for $\widehat{m}_{\pm}(x)$, or $\delta = o(h^{-1/2})$, for $\widetilde{m}_a(x)$, although it is argued in Remark 5 that the new estimators have good asymptotic properties, there can be adverse effects on the finite sample biases. In applications caution is needed when employing constant $\delta$. Each of $\widetilde{m}_{\pm}(x)$ and $\widetilde{m}_a(x)$ uses observations further away from $x$ when using larger values of $\delta$ and that may substantially increase the finite sample bias. In general, larger values of $\delta$, for example, $\delta \in [1.5, 2]$, are recommended only when $n$ is large or when the estimation problem is not difficult, that is, the regression function is smooth and the noise level is low. Otherwise, a smaller $\delta$ is preferred. Furthermore, $\delta = 1$ is a good default. Results of a numerical study, reported in Section 5, support these suggestions.

4.3. *Computation.* A naive way to compute the proposed estimators is to calculate the required local linear estimators and then form the linear combinations. Then the computational effort is increased by some constant factors. This extra burden can be avoided by a careful implementation. Suppose the estimators are computed over an equispaced grid. Letting the spacing of the grid be a multiple of $\omega_n$ reduces the number of kernel evaluations to essentially the same as required by $\widehat{m}$. In addition, fast implementation of kernel estimators, for example, fast Fourier transform or binning methods of Fan and Marron [12], can be employed to alleviate the computational effort.

**5. Numerical performance.** A simulation study was carried out to investigate the finite sample performances. Three regression functions (see [24]),

1. bimodal, $m(x) = 0.3\exp\{-16(x-0.25)^2\} + 0.7\exp\{-64(x-0.75)^2\}$,
2. linear with peak, $m(x) = 2 - 5x + 5\exp\{-400(x-0.5)^2\}$,
3. sine, $m(x) = \sin(5\pi x)$,

were considered. The design was Uniform$(0,1)$. The random error $\varepsilon$ was Normal$(0,1)$ distributed and $\sigma(x) = k\sigma_0$ for all $x \in [0,1]$, where $k = 0.5, 1$ or 2 and $\sigma_0 = 0.1, \sqrt{0.5}$ and 0.5 for the bimodal, linear with peak and sine regression functions, respectively. The sample size was $25, 50, 100, 250$ or $500$. To avoid the data sparseness problem, let $\widehat{m}_{\mathrm{HT}}(x)$ be the local linear estimator employing the interpolation method of Hall and Turlach [15] with the parameter $r$ therein being taken as 3. Let $\widehat{m}_{\mathrm{CH}}(x)$ denote the bias reduction method of Choi and Hall [7] applied to $\widehat{m}_{\mathrm{HT}}(x)$ rather than the local linear estimator. Our variance reduction methods, using $\delta = 0.6$, $0.8$, $1$, $1.2$ or $1.6$, were applied to $\widehat{m}_{\mathrm{HT}}(x)$. The Epanechnikov kernel was employed. The bandwidth $h$ ranged over $\{0.008 \cdot 1.1^k, k = 0, 1, \ldots, 40\}$. In each setting 1000 samples were simulated. The mean integrated squared error (MISE) of an estimator is approximated by the average of the 1000 numerical integrals of the squared errors. The integrated squared bias (ISB) and integrated variance (IV) are likewise approximated.

Figure 3 depicts the MISE, ISB and IV curves under one of the configurations. Our estimator $\widetilde{m}_a$ significantly improves on $\widehat{m}_{\mathrm{HT}}$ in terms of IV. The second-order bias effect of $\widetilde{m}_a$ is more apparent when $h$ is large or $\delta$ is large. The ISB curve of $\widehat{m}_{\mathrm{CH}}$ is much smaller than that of $\widehat{m}_{\mathrm{HT}}$ except when $h$ is large.

Finite sample efficiency of an estimator relative to $\widehat{m}_{\mathrm{HT}}$ can be measured by the ratio of the minimal, over $h$, MISE values. Figure 4 plots some of the results. Our estimator $\widetilde{m}_a$, using any of the $\delta$ values, outperforms $\widehat{m}_{\mathrm{HT}}$ most of the time. The bias reduction estimator $\widehat{m}_{\mathrm{CH}}$ is better than $\widetilde{m}_a$ with $\delta \leq 1$ when the regression function is smooth and when $n \geq 100$. Interestingly, even when $n = 500$, $\widehat{m}_{\mathrm{CH}}$ has no significant advantage over $\widetilde{m}_a$ with $\delta = 1.2$ or $1.6$. Although we do not report the simulation results on the linear with peak regression function, we observe that, in this case, $\widehat{m}_{\mathrm{CH}}$ performs roughly the same as $\widehat{m}_{\mathrm{HT}}$, but our variance reduction estimator $\widetilde{m}_a$ outperforms $\widehat{m}_{\mathrm{HT}}$. These observations coincide with previous studies on bias reduction methods that they usually depend on certain higher-order approximations which take effect only when the curve is smooth and $n$ is large. In the high-noise case ($k = 2$), $\widehat{m}_{\mathrm{HT}}$ breaks down and there is not much improvement employing any modification.

Examining Figure 4 more closely, we have the following conclusions. First, $\delta = 1$ is a reliable choice for general purposes. Except when the noise level is high ($k = 2$), the relative efficiency for $\delta = 1$ is already above 1.1 from $n = 100$ and is not far from the asymptotic value $\gamma_a(1) \approx 1.22$ (see Figure 1) for moderate sample sizes. In all of the cases, the curves for $\delta = 0.6$ and $0.8$ become almost flat starting from $n = 100$. Besides $\delta = 1$, larger $\delta$ values, for example, $\delta = 1.2$, may have potential in practice. In general, for smooth curves like the sine function or low noise levels, gains resulted from employing $\delta > 1$ over using $\delta \leq 1$ become noticeable for $n \geq 500$.

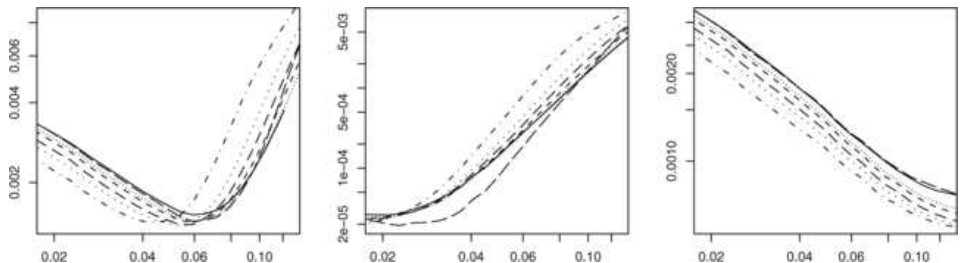

FIG. 3. *MISE, ISB and IV curves. From left to right are MISE, ISB and IV plotted against $h$ when the regression is bimodal, the design is $\mathrm{U}(0,1)$, $k = 1$ and $n = 100$. The line types solid, long-dashed, dotted, short-dashed, dashed, dotted-spaced and dotted-dashed respectively represent $\widehat{m}_{\mathrm{HT}}$, $\widehat{m}_{\mathrm{CH}}$ and $\widetilde{m}_a$ with $\delta = 0.6, 0.8, 1, 1.2$ and $1.6$.*

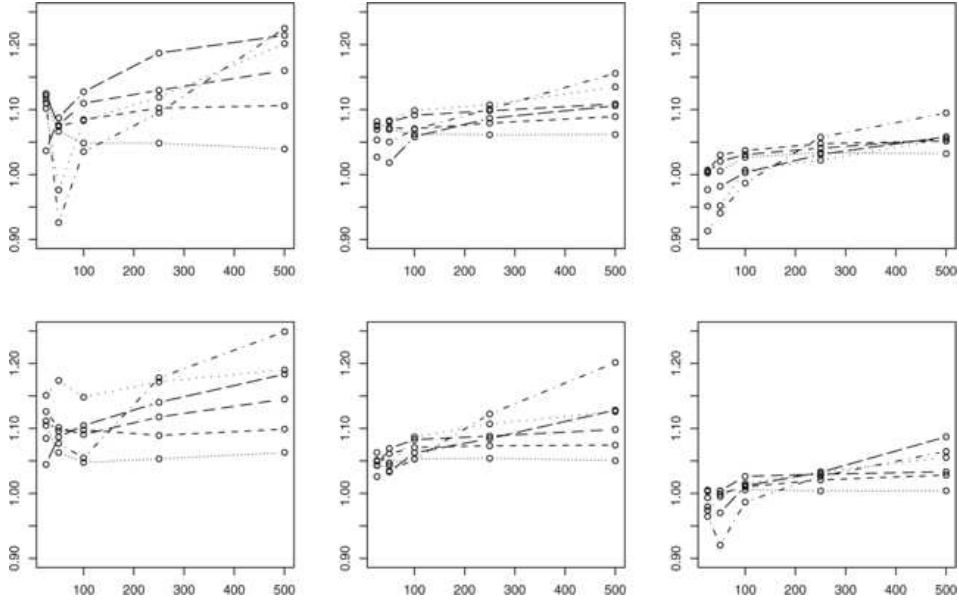

Fig. 4. *MISE relative efficiency. Top (or bottom) row plots against n the MISE efficiencies relative to $\widehat{m}_{HT}$ when regression is the bimodal (or sine) function and design is $U(0,1)$. From left, the columns correspond to the noise levels $k = 0.5, 1$ and $2$. Line types represent the estimators and are as in Figure 3.*

Although the simulation results of $\widetilde{m}_{\pm}$ are not given here, the performances of $\widetilde{m}_{\pm}$ are hard to differentiate from each other, the average version $\widetilde{m}_a$ copes with the second-order bias effect better than $\widetilde{m}_{\pm}$, and $\widetilde{m}_{\pm}(x)$ using $\delta \leq 1.2$ outperform $\widehat{m}_{HT}$ most of the time. The simulation study has been summarized in terms of the MISE performances. The pointwise mean squared errors, not reported here, provide similar conclusions. In addition, the MSE relative efficiencies of the proposed estimators are roughly unchanged as $x$ varies. Two truncated Normal designs, $N(0.5, 0.5^2) \cap (0,1)$ and $N(0,1) \cap (0,1)$, were experimented with under all the considered configurations as well. In all cases, the relative MISE efficiency of each of our variance reduction estimators behaves consistently across the three different designs.

**6. Generalizations and applications.** Generalizations of the methodology are multifold and comparatively easy. First, in kernel estimation of curves and their derivatives using higher-order kernels or by fitting higher-degree local polynomials, a linear combination of some preliminary estimators at more than three neighboring points is required in order that moment conditions, in the same spirit as those in (7.4), are satisfied. See also the asymptotic mean expressions (7.3) and (7.5). Second, extension to estimation of multi-dimensional surfaces is possible. Cheng and Peng [6] made some progress in this regard: when using product kernels, form a grid in every direction as in the one-dimensional case, and then take the coefficient of one preliminary estimator as the product of all the corresponding one-dimensional coefficients. Asymptotic relative efficiency, in terms of MSE or MISE, compared to the local linear smoother is $\gamma_q(\delta)^d$ and is $\gamma_a(\delta)^d$ for the average version.

The proposed variance reduction strategy is fairly simple, allowing a wide range of potential applications. Variance reduction is particularly useful when very few data points can be used in the estimation, for example, estimating quantities in conditional distributions. Other useful applications include various local modeling techniques such as local likelihood estimation, varying coefficient models and hazard regression. In situations where the covariance structure of the preliminary estimators is completely analogous to (7.6), direct application is appropriate. Examples include kernel density estimation and are found in local likelihood modeling [19, 26]. Notice that applying the techniques to density estimation can introduce negativity. When applying to more complicated settings, the covariance structure of the preliminary estimators needs to be investigated.

**7. Proofs.**

PROOF OF THEOREM 1. Following [9], consider the version of the local linear estimator in (3.1). Denote $Z_n = O_l(a_n)$ if $E|Z_n|^l = O(a_n^l)$ and $Z_n = o_l(a_n)$ if $E|Z_n|^l = o(a_n^l)$. For any fixed $z \in [\alpha_{x,0}, \alpha_{x,2}]$, write $\omega_{i,z} = hK_h(z - X_i)\{S_{n,2}(z) - (z - X_i)S_{n,1}(z)\}$ with $S_{n,j}(z) = h\sum_{i=1}^{n}(z - X_i)^j K_h(z - X_i)$, $j = 0, 1, 2$. Then one can show that

$$(7.1) \quad n^2 h^4 \left\{\sum_{i=1}^{n}\omega_{i,z} + n^{-2}\right\}^{-1} = \{\nu_{20}f^2(z)\}^{-1} + o_4(1),$$

$$(7.2) \quad \sum_{i=1}^{n}\{m(X_i) - m(z)\}\omega_{i,z} = n^2 h^6 f^2(z)\nu_{20}S_{n,z} + o_4(n^2 h^6),$$

where $S_{n,z} = h^{-2}\mathrm{E}\{m(X) - m(z) - m'(z)(X - z)\}K_h(z - X)$. From (7.1) and (7.2),

$$(7.3) \quad \mathrm{E}\{\widehat{m}(z)\} = m(z) + \tfrac{1}{2}m''(z)\nu_{20}h^2 + o\{h^2 + (nh)^{-1/2}\}.$$

Note that, for $j = 0, 1, 2$,

(7.4) $$A_0(r)(-1-r)^j + A_1(r)(-r)^j + A_2(r)(1-r)^j = \delta_{0,j},$$

where $\delta_{0,j} = 1$ if $j = 0$ and 0 otherwise. Recall that $x = \alpha_{x,1} + r\omega_n \in [\alpha_{x,0}, \alpha_{x,2}]$, $-1 \leq r \leq 1$, and $\omega_n = \delta h = \alpha_{x,2} - \alpha_{x,1} = \alpha_{x,1} - \alpha_{x,0}$. Using (7.3), Taylor expansion and (7.4) we have

$$\mathrm{E}\{\widetilde{m}_q(x)\} = \sum_{i=0,1,2} A_i(r) \mathrm{E}\{\widehat{m}(\alpha_{x,i})\}$$

$$= \{m(x) + \tfrac{1}{2} m''(x)\nu_{20} h^2\} \sum_{i=0,1,2} A_i(r)$$

(7.5) $$+ m'(x) \sum_{i=0,1,2} A_i(r)(\alpha_{x,i} - x)$$

$$+ \tfrac{1}{2} m''(x) \sum_{i=0,1,2} A_i(r)(\alpha_{x,i} - x)^2 + o\{h^2 + (nh)^{-1/2}\}$$

$$= m(x) + \tfrac{1}{2} m''(x)\nu_{20} h^2 + o\{h^2 + (nh)^{-1/2}\}.$$

Next we compute the variance of $\widetilde{m}_q(x)$. For any $u, v \in [\alpha_{x,0}, \alpha_{x,2}]$, from (7.1) and (7.3),

$$\mathrm{Cov}\{\widehat{m}(u), \widehat{m}(v)\} = \mathrm{Cov}\left\{\frac{\sum_{i=1}^{n} \omega_{i,u}\{Y_i - m(u)\}}{\sum_{i=1}^{n} \omega_{i,u} + n^{-2}}, \frac{\sum_{i=1}^{n} \omega_{i,v}\{Y_i - m(v)\}}{\sum_{i=1}^{n} \omega_{i,v} + n^{-2}}\right\}$$

$$+ o\left(\frac{1}{n^4 h^4}\right).$$

Taking conditional expectations on $X_1, \ldots, X_n$ and using mean and variance decomposition yield

$$\mathrm{Cov}\{\widehat{m}(u), \widehat{m}(v)\}$$

$$= \mathrm{E}\left[\frac{\sum_{i=1}^{n} \omega_{i,u}\{m(X_i) - m(u)\} \sum_{j=1}^{n} \omega_{j,v}\{m(X_j) - m(v)\}}{(\sum_{i=1}^{n} \omega_{i,u} + n^{-2})(\sum_{j=1}^{n} \omega_{j,v} + n^{-2})}\right]$$

(7.6) $$- \mathrm{E}\left[\frac{\sum_{i=1}^{n} \omega_{i,u}\{m(X_i) - m(u)\}}{\sum_{i=1}^{n} \omega_{i,u} + n^{-2}}\right] \mathrm{E}\left[\frac{\sum_{j=1}^{n} \omega_{j,v}\{m(X_j) - m(v)\}}{\sum_{j=1}^{n} \omega_{j,v} + n^{-2}}\right]$$

$$+ \mathrm{E}\left[\frac{\sum_{i=1}^{n} \omega_{i,u}\omega_{i,v}\sigma^2(X_i)}{(\sum_{i=1}^{n} \omega_{i,u} + n^{-2})(\sum_{i=1}^{n} \omega_{i,v} + n^{-2})}\right] + o\{(nh)^{-4}\}$$

$$= \frac{\sigma^2(x)}{nf(x)} \int K_h(u-t) K_h(v-t)\, dt + o\{h^4 + (nh)^{-1}\},$$

where the last equality follows from (7.1), (7.2) and $\sum_{i=1}^{n} \omega_{i,u}\omega_{i,v}\sigma^2(X_i) = n^3 h^8 \sigma^2(x) \nu_{20}^2 f^2(x) \int K_h(u-t) K_h(v-t) f(t)\, dt \{1 + O_l(h^\alpha + (nh)^{-1/2})\}$.

Then (3.5) is valid since

$$\operatorname{Var}\{\widetilde{m}_q(x)\} = \sum_{i=0,1,2} A_i(r)^2 \operatorname{Var}\{\widehat{m}(\alpha_{x,i})\}$$
$$+ \sum_{i=0,1,2} \sum_{j \neq i} A_i(r) A_j(r) \operatorname{Cov}\{\widehat{m}(\alpha_{x,i}), \widehat{m}(\alpha_{x,j})\}. \qquad \square$$

PROOF OF PROPOSITION 1. Property (a) follows from $\int K(x-\delta/2) K(x+\delta/2)\,dx = \int K(x) K(x-\delta)\,dx = \int K(x) K(x+\delta)\,dx$ and writing

$$C(\delta) = \int \{\tfrac{3}{2} K(x)^2 - K(x) K(x+\delta) - K(x) K(x-\delta)$$
(7.7)
$$+ \tfrac{1}{2} K(x-\delta) K(x+\delta)\}\,dx$$
$$= \int \{K(x) - \tfrac{1}{2} K(x+\delta) - \tfrac{1}{2} K(x-\delta)\}^2\,dx.$$

Property (b) can be shown by differentiating the right-hand side of (7.7). $\square$

PROOF OF PROPOSITION 2. From (3.7) and (3.9),

$$D(\delta) = \frac{2nhf(x)}{\sigma^2(x)} \left[\frac{1}{2} \operatorname{Var}\{\widetilde{m}_+(x)\} + \frac{1}{2} \operatorname{Var}\{\widetilde{m}_-(x)\} - \operatorname{Var}\{\widetilde{m}_a(x)\}\right]$$
$$\times \{1 + o(1)\}$$
$$= \frac{nhf(x)}{2\sigma^2(x)} \operatorname{Var}\{\widetilde{m}_+(x) - \widetilde{m}_-(x)\}\{1 + o(1)\}. \qquad \square$$

PROOF OF THEOREM 3. Let $\mu_{ijk}(x) = \operatorname{E}\{w_{ijk}(x)\}$, $\Delta_{ijk}(x) = w_{ijk}(x) - \mu_{ijk}(x)$, $\mu^*(x) = \mu_{210}(x)\mu_{011}(x) - \mu_{110}(x)\mu_{111}(x)$, $U_1(x) = w_{210}(x)w_{011}(x) - w_{110}(x)w_{111}(x)$ and $U_2(x) = w_{210}(x)w_{010}(x) - w_{110}^2(x)$. Note that

$$U_1(x) = \Delta_{210}(x)\Delta_{011}(x) - \Delta_{110}(x)\Delta_{111}(x) + \mu_{011}(x)\Delta_{210}(x)$$
$$+ \mu_{210}(x)\Delta_{011}(x) - \mu_{110}(x)\Delta_{111}(x) - \mu_{111}(x)\Delta_{110}(x) + \mu^*(x),$$
$$U_2(x) = \Delta_{210}(x)\Delta_{010}(x) - \Delta_{110}^2(x) - 2\mu_{110}(x)\Delta_{110}(x) - \mu_{110}^2(x)$$
$$+ \mu_{210}(x)\Delta_{010}(x) + \mu_{010}(x)\Delta_{210}(x) + \mu_{210}(x)\mu_{010}(x),$$

where $\Delta_{210}(x)$, $\Delta_{011}(x)$, $\Delta_{110}(x)$, $\Delta_{111}(x)$, $\Delta_{010}(x)$ and $\Delta_{012}(x)$ are all $O_p\{(nh)^{-1/2}\}$ and $\mu_{210}(x) = O(1)$, $\mu_{011}(x) = O(h^2)$, $\mu_{110}(x) = O(h)$, $\mu_{111}(x) = O(h)$, $\mu_{010}(x) = O(1)$ and $\mu_{012}(x) = O(1)$. Then

$$T_n(x) = \frac{\sqrt{nh}}{\sqrt{\nu_{02}}} \frac{U_1(x)}{U_2(x)} w_{010}(x) \left\{ w_{012}(x) - 2w_{011}(x)\frac{U_1(x)}{U_2(x)} + w_{010}(x)\frac{U_1^2(x)}{U_2^2(x)} \right\}^{-1/2}$$

$$= \frac{\sqrt{nh}}{\sqrt{\nu_{02}}} U_1(x) w_{010}(x)$$

$$\times \{w_{012}(x)U_2^2(x) - 2w_{011}(x)U_1(x)U_2(x) + w_{010}(x)U_1^2(x)\}^{-1/2}$$

$$= \frac{\sqrt{nh}}{\sqrt{\nu_{02}}} \{U_1(x)\Delta_{010}(x) + \mu_{010}(x)U_1(x)\}$$

$$\times \{U_2^2(x)\Delta_{012}(x) + \mu_{012}(x)U_2^2(x)$$
$$- 2\Delta_{011}(x)U_1(x)U_2(x) - 2\mu_{011}(x)U_1(x)U_2(x)$$
$$+ \Delta_{010}(x)U_1^2(x) + \mu_{010}(x)U_1^2(x)\}^{-1/2}$$

$$= \frac{\sqrt{nh}}{\sqrt{\nu_{02}}} [\mu_{210}(x)\Delta_{011}(x)\Delta_{010}(x) + \mu_{010}(x)\Delta_{210}(x)\Delta_{011}(x)$$
$$- \mu_{010}(x)\Delta_{110}(x)\Delta_{111}(x) + \mu_{010}(x)\mu_{210}(x)\Delta_{011}(x)$$
$$- \mu_{010}(x)\mu_{110}(x)\Delta_{111}(x) - \mu_{010}(x)\mu_{111}(x)\Delta_{110}(x)$$
$$+ \mu_{010}(x)\mu^*(x) + O_p(\{(nh)^{-1/2} + h\}^3)]$$
$$\times [\mu_{210}^2(x)\mu_{010}^2(x)\Delta_{012}(x) + \mu_{210}^2(x)\mu_{010}^2(x)\mu_{012}(x)$$
$$+ 2\mu_{210}^2(x)\mu_{010}(x)\mu_{012}(x)\Delta_{010}(x)$$
$$+ 2\mu_{210}(x)\mu_{010}^2(x)\mu_{012}(x)\Delta_{210}(x)$$
$$+ O_p(\{(nh)^{-1/2} + h\}^2)]^{-1/2}$$

$$= T_{n1}(x) + O_p(\{(nh)^{-1/2} + h\}^2),$$

where

$$T_{n1}(x) = \{\mu_{210}^2(x)\mu_{010}^2(x)\mu_{012}(x)\}^{-1/2}\nu_{02}^{-1/2}$$
$$\times \sqrt{nh}\{\mu_{210}(x)\Delta_{011}(x)\Delta_{010}(x) + \mu_{010}(x)\Delta_{210}(x)\Delta_{011}(x)$$
$$- \mu_{010}(x)\Delta_{110}(x)\Delta_{111}(x) + \mu_{010}(x)\mu_{210}(x)\Delta_{011}(x)$$
$$- \mu_{010}(x)\mu_{110}(x)\Delta_{111}(x) - \mu_{010}(x)\mu_{111}(x)\Delta_{110}(x)$$
$$+ \mu_{010}(x)\mu^*(x)\}$$
$$- 2^{-1}\{\mu_{210}^2(x)\mu_{010}^2(x)\mu_{012}(x)\}^{-3/2}\nu_{02}^{-1/2}\sqrt{nh}\mu_{210}^2(x)\mu_{010}^2(x)$$
$$\times \{\mu_{210}(x)\mu_{010}(x)\Delta_{011}(x)\Delta_{012}(x)$$
$$+ 2\mu_{210}(x)\mu_{012}(x)\Delta_{011}(x)\Delta_{010}(x)$$
$$+ 2\mu_{010}(x)\mu_{012}(x)\Delta_{011}(x)\Delta_{210}(x)\}.$$

We have
$$\mathrm{E}\{T_{n1}(x)\} = \sqrt{nh}\{\mu_{210}^2(x)\mu_{010}^2(x)\mu_{012}(x)\}^{-1/2}\nu_{02}^{-1/2}\mu_{010}(x)\mu^*(x)$$

$$- 2^{-1}(nh)^{-1/2}\{\mu_{210}^2(x)\mu_{010}^2(x)\mu_{012}(x)\}^{-3/2}$$
$$\times \nu_{02}^{-1/2}\mu_{210}^3(x)\mu_{010}^3(x)\mu_{023}(x)$$
$$+ O(\{(nh)^{-1/2} + h\}^2)$$
$$= 2^{-1}\nu_{21}f^{1/2}(x)\nu_{02}^{-1/2}\sigma^{-1}(x)m''(x)\sqrt{nh}h^2$$
$$- 2^{-1}(nh)^{-1/2}f^{-1/2}(x)\nu_{02}^{1/2}\sigma^{-3}(x)V_3(x) + O(\{(nh)^{-1/2} + h\}^2),$$

$$\mathrm{E}\{T_{n1}^2(x)\} = 1 + O(\{(nh)^{-1/2} + h\}^2),$$

$$\mathrm{E}\{T_{n1}^3(x)\} = (nh)^{-1/2}\{\mu_{210}^2(x)\mu_{010}^2(x)\mu_{012}(x)\nu_{02}\}^{-3/2}\mu_{010}^3(x)\mu_{210}^3(x)$$
$$\times \{\mu_{033}(x) - \tfrac{9}{2}\mu_{022}(x)\mu_{023}(x)/\mu_{012}(x)\} + O(\{(nh)^{-1/2} + h\}^2)$$
$$= (nh)^{-1/2}f^{-1/2}(x)\sigma^{-3}(x)\nu_{02}^{-3/2}V_3(x)\{\nu_{03} - \tfrac{9}{2}\nu_{02}^2\}$$
$$+ O(\{(nh)^{-1/2} + h\}^2),$$

$$\mathrm{E}\{T_{n1}^l(x)\} = O(\{(nh)^{-1/2} + h\}^2) \quad \text{for } l \geq 4.$$

Hence, by Edgeworth expansion (see, e.g., Chapter 2 of [14]),

$$P\{T_{n1}(x) \leq z\}$$
$$= \Phi(z) + (nh^5)^{1/2}4^{-1}\nu_{21}\nu_{02}^{-1/2}\sigma^{-1}(x)f^{1/2}(x)m''(x)(z^2 - 3)\phi(z)$$
$$- (nh)^{-1/2}6^{-1}\nu_{02}^{-3/2}\sigma^{-3}(x)f^{-1/2}(x)V_3(x)\{\nu_{03}(z^2 - 1) - 3\nu_{02}^2 z^2\}\phi(z)$$
$$+ O(\{(nh)^{-1/2} + h\}^2),$$

and then applying the delta method yields (3.14). To calculate the coverage probability of $\widetilde{I}_\beta$, write

$$T_n^*(x) = \frac{\sqrt{nh}}{\sqrt{\widetilde{\nu}_{02}}}\left\{\sum_{l=0,1,2} A_l(r)\frac{\widehat{m}(\alpha_{x,l}) - m(\alpha_{x,l})}{\sqrt{\widehat{\sigma}^2(\alpha_{x,l})/w_{010}^2(\alpha_{x,l})}}\frac{\sqrt{\widehat{\sigma}^2(\alpha_{x,l})/w_{010}^2(\alpha_{x,l})}}{\sqrt{\widehat{\sigma}^2(x)/w_{010}^2(x)}}\right.$$
$$\left. + \sum_{l=0,1,2} A_l(r)\frac{m(\alpha_{x,l}) - m(x)}{\sqrt{\widehat{\sigma}^2(x)/w_{010}^2(x)}}\right\}$$
$$= \sum_{l=0,1,2} A_l(r)\left[\sqrt{\nu_{02}/\widetilde{\nu}_{02}}T_n(\alpha_{x,l})\right.$$
$$\left. + \sqrt{\nu_{02}/\widetilde{\nu}_{02}}T_n(\alpha_{x,l})\left\{\frac{\sqrt{\widehat{\sigma}^2(\alpha_{x,l})/w_{010}^2(\alpha_{x,l})}}{\sqrt{\widehat{\sigma}^2(x)/w_{010}^2(x)}} - 1\right\}\right.$$

$$+ \frac{m(\alpha_{x,l}) - m(x)}{\sqrt{\hat{\sigma}^2(x)/w_{010}^2(x)}}\Bigg].$$

Since
$$T_n(\alpha_{x,l}) = T_{n1}(\alpha_{x,l}) + O_p(\{(nh)^{-1/2} + h\}^2),$$
$$T_{n1}(\alpha_{x,l}) = \{\mu_{210}^2(\alpha_{x,l})\mu_{010}^2(\alpha_{x,l})\mu_{012}(\alpha_{x,l})\}^{-1/2}$$
$$\times \sqrt{nh}\nu_{02}^{-1/2}\mu_{010}(\alpha_{x,l})\mu_{210}(\alpha_{x,l})\Delta_{011}(\alpha_{x,l})$$
$$\times [1 + O_p\{(nh)^{-1/2} + h\}],$$

we have
$$T_n^*(x) = T_{n1}^*(x) + O_p(\{(nh)^{-1/2} + h\}^2),$$

where
$$T_{n1}^*(x) = \sum_{l=0,1,2} A_l(t)[\sqrt{\nu_{02}/\widetilde{\nu}_{02}}T_{n1}(\alpha_{x,l})$$
$$- 2^{-1}\mu_{012}^{-1}(x)\mu_{010}^2(x)$$
$$\times \{\mu_{012}(\alpha_{x,l})\mu_{010}^{-2}(\alpha_{x,l}) - \mu_{012}(x)\mu_{010}^{-2}(x)\}$$
$$\times \{\mu_{210}^2(\alpha_{x,l})\mu_{010}^2(\alpha_{x,l})\mu_{012}(\alpha_{x,l})\}^{-1/2}$$
$$\times \sqrt{nh}\widetilde{\nu}_{02}^{-1/2}\mu_{010}(\alpha_{x,l})\mu_{010}(\alpha_{x,l})\Delta_{011}(\alpha_{x,l})].$$

Then the following equalities and the delta method yield (3.15).
$$E\{T_{n1}^*(x)\} = 2^{-1}\nu_{21}f(x)^{1/2}\widetilde{\nu}_{02}^{-1/2}\sigma^{-1}(x)m''(x)\sqrt{nh^5}$$
$$- 2^{-1}f(x)^{-1/2}\widetilde{\nu}_{02}^{1/2}\sigma^{-3}(x)V_3(x)(nh)^{-1/2}$$
$$+ O(\{(nh)^{-1/2} + h\}^2),$$
$$E\{T_{n1}^*(x)^2\} = 1 + O(\{(nh)^{-1/2} + h\}^2),$$
$$E\{T_{n1}^*(x)^3\} = (nh)^{-1/2}f(x)^{-1/2}\sigma^{-3}(x)\widetilde{\nu}_{02}^{-3/2}V_3(x)\{\widetilde{\nu}_{03} - 2^{-1}9\widetilde{\nu}_{02}^2\}$$
$$+ O(\{(nh)^{-1/2} + h\}^2),$$
$$E\{T_{n1}^*(x)^l\} = O(\{(nh)^{-1/2} + h\}^2) \quad \text{for } l \geq 4. \qquad \square$$

**Acknowledgments.** The authors would like to thank one of the Co-Editors, an Associate Editor and the referees for constructive comments that led to substantial improvements of the paper. Mr. Lu-Hung Chen's help with the numerical work is acknowledged.

M.-Y. CHENG
DEPARTMENT OF MATHEMATICS
NATIONAL TAIWAN UNIVERSITY
TAIPEI 106
TAIWAN
E-MAIL: cheng@math.ntu.edu.tw
URL: www.math.ntu.edu.tw/cheng/

L. PENG
SCHOOL OF MATHEMATICS
GEORGIA INSTITUTE OF TECHNOLOGY
ATLANTA, GEORGIA 30332
USA
E-MAIL: peng@math.gatech.edu
URL: www.math.gatech.edu/~peng/

J.-S. WU
DEPARTMENT OF MATHEMATICS
TAMKANG UNIVERSITY
TAMSUI 251
TAIWAN
E-MAIL: jswu@math.tku.edu.tw